\documentclass[12pt]{amsart}
\begin{document}
\numberwithin{equation}{section}

\def\1#1{\overline{#1}}
\def\2#1{\widetilde{#1}}
\def\3#1{\widehat{#1}}
\def\4#1{\mathbb{#1}}
\def\5#1{\frak{#1}}
\def\6#1{{\mathcal{#1}}}

\def\C{{\4C}}
\def\R{{\4R}}
\def\N{{\4N}}
\def\Z{{\4Z}}
\def\K{{\4K}}

\title[Formal and finite order equivalences]{Formal and finite order equivalences}
\author[D. Zaitsev]{Dmitri Zaitsev*}
\thanks{*Supported in part by the Science Foundation Ireland grant 06/RFP/MAT018.}
\address{D. Zaitsev: School of Mathematics, Trinity College Dublin, Dublin 2, Ireland}
\email{zaitsev@maths.tcd.ie}

\begin{abstract} 
We show that two families of germs of real-analytic subsets in $\C^{n}$
are formally equivalent if and only if they are equivalent of any finite order.
We further apply the same technique to obtain analogous statements
for equivalences of real-analytic self-maps and vector fields under conjugations.
On the other hand, we provide an example of two sets of germs of smooth curves 
that are equivalent of any finite order but not formally equivalent.
\end{abstract}

\maketitle

\def\Label#1{\label{#1}}


\def\cn{{\C^n}}
\def\cnn{{\C^{n'}}}
\def\ocn{\2{\C^n}}
\def\ocnn{\2{\C^{n'}}}


\def\dist{{\rm dist}}
\def\const{{\rm const}}
\def\rk{{\rm rank\,}}
\def\id{{\sf id}}
\def\aut{{\sf aut}}
\def\Aut{{\sf Aut}}
\def\CR{{\rm CR}}
\def\GL{{\sf GL}}
\def\Re{{\sf Re}\,}
\def\Im{{\sf Im}\,}
\def\span{\text{\rm span}}

\def\codim{{\rm codim}}
\def\crd{\dim_{{\rm CR}}}
\def\crc{{\rm codim_{CR}}}

\def\phi{\varphi}
\def\eps{\varepsilon}
\def\d{\partial}
\def\a{\alpha}
\def\b{\beta}
\def\g{\gamma}
\def\G{\Gamma}
\def\D{\Delta}
\def\Om{\Omega}
\def\k{\kappa}
\def\l{\lambda}
\def\L{\Lambda}
\def\z{{\bar z}}
\def\w{{\bar w}}
\def\Z{{\mathbb Z}}
\def\t{\tau}
\def\th{\theta}

\emergencystretch15pt
\frenchspacing

\newtheorem{Thm}{Theorem}[section]
\newtheorem{Cor}[Thm]{Corollary}
\newtheorem{Pro}[Thm]{Proposition}
\newtheorem{Lem}[Thm]{Lemma}

\theoremstyle{definition}\newtheorem{Def}[Thm]{Definition}

\theoremstyle{remark}
\newtheorem{Rem}[Thm]{Remark}
\newtheorem{Exa}[Thm]{Example}
\newtheorem{Exs}[Thm]{Examples}

\def\bl{\begin{Lem}}
\def\el{\end{Lem}}
\def\bp{\begin{Pro}}
\def\ep{\end{Pro}}
\def\bt{\begin{Thm}}
\def\et{\end{Thm}}
\def\bc{\begin{Cor}}
\def\ec{\end{Cor}}
\def\bd{\begin{Def}}
\def\ed{\end{Def}}
\def\br{\begin{Rem}}
\def\er{\end{Rem}}
\def\be{\begin{Exa}}
\def\ee{\end{Exa}}
\def\bpf{\begin{proof}}
\def\epf{\end{proof}}
\def\ben{\begin{enumerate}}
\def\een{\end{enumerate}}
\def\beq{\begin{equation}}
\def\eeq{\end{equation}}

\section{Introduction}

There are three basic equivalence relations between germs of real-analytic submanifolds in $\C^{n}$. The first is that of {\em biholomorphic equivalence}, the second of {\em formal equivalence} and the third of {\em equivalence of any finite order} (see \S\ref{preli} below for precise definitions). In \cite{BRZ} Baouendi, Rothschild and the author proved that at {\em points of general position}, these three notions of equivalence coincide
(where an equivalence of order $k$ was called a $k$-equivalence). Moser and Webster~\cite{MW} gave an example of two germs of surfaces in $\C^{2}$ (at their complex points) that are formally but not biholomorphically equivalent.
Thus, in general, biholomorphic and formal equivalences do not coincide.
It remained a question whether formal equivalence for real-analytic submanifolds always coincides with 
their equivalence of any finite order.

The present paper answers the latter question affirmatively. 
In fact, the affirmative answer is given even for (possibly singular)
real-analytic subsets as well as for their arbitrary families (Theorem~\ref{maincor}).
On the other hand, in \S\ref{example} we provide an example showing
that for countable unions of smooth complex curves in $\C^{2}$,
the notions of formal equivalence and equivalence of any finite order do not coincide in general.

The technique used in the proof can be applied in other situations.
In \S\ref{fields} we demonstrate it for equivalence relations between real-analytic
(in particular, also holomorphic) self-maps under biholomorphic conjugations
as well as for closely related equivalence relations between vector fields.

\section{Preliminaries}\label{preli}
In the following we summarize the definitions for various equivalence relations.
Let $\Phi\colon (\C^n,0)\to (\C^n,0)$ be a formal invertible map given by 
a power series in $z=(z_{1},\ldots,z_{n})\in\C^{n}$
vanishing at $0$.

\subsection{Equivalences for ideals and real-analytic sets}
We write $z=x+iy\in \R^{n}+i\R^{n}$ and denote by $\R[[x,y]]$
the ring of all formal power series in $x$ and $y$.
As customary, denote by $\5m^{k}\subset \R[[x,y]]$
the maximal ideal consisting or all power series vanishing at $0$.

\begin{enumerate}
\item
$\Phi$ is said to be a {\em formal equivalence} between two subsets
$I,J\subset \R[[x,y]]$ if $f\circ\Phi\in I$ for every $f\in J$ and 
$ g\circ (\Phi^{-1})\in J$ for every $g\in I$.
\item $\Phi$ is said to be an {\em equivalence of order $k$}  between two subets
$I,J\subset \R[[x,y]]$
if $f\circ \Phi\in I+\5m^k$ for every $f\in J$ and 
$ g\circ (\Phi^{-1})+\5m^k\in J$ for every $g\in I$.
\item $\Phi$ is a said to be a {\em formal equivalence} (resp.\ {\em equivalence of order $k$}) between two germs at $0$ of real-analytic
subsets $S,T\subset\C^n$ if it is a formal equivalence (resp.\ equivalence of order $k$) between
their ideals in $\R[[x,y]]$. We write $\Phi_*S=T$ (resp.\ $\Phi_*S\sim_k T$).
\item $\Phi$ is said to be a {\em formal equivalence} (resp.\ {\em equivalence of order $k$}) between two families
$(S_\a)_{\a\in A}$ and $(T_\a)_{\a\in A}$ of germs at $0$ of real-analytic
subsets in $\C^n$ if $\Phi_* S_\a = T_\a$  (resp.\ $\Phi_* S_\a \sim_k T_\a$)  for all $\a\in A$.
\item $\Phi$ is said to be a {\em formal equivalence} (resp.\ {\em equivalence of order $k$}) between two sets
$\{S_\a\}_{\a\in A}$ and $\{T_\b\}_{\b\in B}$ of germs at $0$ of real-analytic
subsets in $\C^n$ if for every $\a\in A$ there exists $\b\in B$ with $\Phi_*S_\a=T_\b$ (resp.\   $\Phi_*S_\a\sim_k T_\b$) and for every $\b\in B$ there exists $\a\in A$ with  $(\Phi^{-1})_*T_\b=S_\a$ (resp.\   $(\Phi^{-1})_*T_\b\sim_k S_\a$).
\end{enumerate}

Similarly to (4) and (5) one defines formal equivalence and equivalence of order $k$
for families $(I_\a)_{\a\in A}$ and $(J_\a)_{\a\in A}$ of ideals in $\R[[x,y]]$ as well for sets of ideals
$\{I_\a\}_{\a\in A}$ and $\{J_\b\}_{\b\in B}$.

We have the following result stating the coincidence of the two notions of equivalence for families of ideals and real-analytic sets:

\bt\Label{maincor}
Two families of ideals in $\R[[x,y]]$ are formally equivalent
if and only if they are equivalent of any finite order.
In particular, two families of germs at $0$ of real-analytic subsets of $\C^{n}$
are formally equivalent
if and only if they are equivalent of any finite order.
\et

On the other hand, the corresponding notions of equivalence do not coincide in general for (countable) sets of real-analytic sets (see \S\ref{example}).

Theorem~\ref{maincor} will be obtained as a direct consequence of the more precise Theorem~\ref{main} below
which gives a description of the set of all formal equivalences between two families of ideals (or real-analytic subsets) as well as their relation with finite order equivalences.

\section{Semi-algebraic sets and Nash groups}\Label{Nash}
The first main ingredient is the theory of semi-algebraic sets and maps
as well as of Nash manifolds and Nash groups.
For the reader's convenience, we recall here the terminology.
For the proofs of the properties of semi-algebraic sets, 
we refer to Benedetti-Risler \cite{BR}.

\begin{Def}\label{def-s-a}
	A subset $V$ of \/ $\R^n$ is called {\em semi-algebraic}
if it admits some representation of the form
$$V = \bigcup_{i=1}^s \bigcap_{j=1}^{r_i} V_{ij}$$
where, for each $i=1,\ldots,s$, and $j=1,\ldots,r_i$,
$V_{ij}$ is either $\{ x\in\R^n : P_{ij}(x)<0 \}$ or
$\{ x\in\R^n : P_{ij}(x)=0 \}$ for a real polynomial $P_{ij}$.
\end{Def}

	As a consequence of the definition it follows that finite unions and
intersections of semi-algebraic sets are always semi-algebraic. Moreover,
closures, boundaries, interiors
and connected components of semi-algebraic sets are semi-algebraic.

\bp\Label{strat}
	Every semi-algebraic set in $\R^n$ admits a 
stratification into a finite disjoint union of semi-algebraic subsets,
each of which is a connected real-analytic submanifold of $\R^n$.
\ep

In particular, a semi-algebraic set has finitely many connected components.
	The natural morphisms in the category of semi-algebraic set are 
semi-algebraic maps:

\begin{Def}\label{map}
	Let $X\subset\R^n$ and $Y\subset\R^n$ be semi-algebraic sets.
A map $f\colon X\to Y$ is called {\em semi-algebraic} if the graph
of $f$ is a semi-algebraic set in $\R^{m+n}$.
\end{Def}

\bt[Tarski-Seidenberg]\Label{TS}
	Let $f\colon X \to Y$ be a semi-algebraic map. Then the image
$f(X)\subset Y$ is a semi-algebraic set.
\et

Combining real-analytic manifolds with semi-algebraic sets
we obtain the categories of Nash manifolds and Nash groups:

\begin{Def}
\begin{enumerate}
\item A {\em Nash map} is a real-analytic function
$f=(f_1,\ldots,f_m)\colon U\to \R^m$ (where $U$ is an open semi-algebraic
subset of $\R^n$) such that for each of the components $f_k$ there is a
nontrivial polynomial $P$ with $P(x_1,\ldots,x_n,f_k(x_1,\ldots,x_n))=0$
for all $(x_1,\ldots,x_n)\in U$.
\item A {\em Nash manifold} $M$ is a real analytic manifold with finitely many
coordinate charts $\phi_i\colon U_i\to V_i$ such that $V_i\subset\R^n$ is
open semi-algebraic for all $i$ and the transition functions are Nash
(a Nash atlas).
\item A {\em Nash group} is a Nash manifold with a group operation
$(x,y)\to xy^{-1}$ which is Nash with respect to every Nash coordinate chart.
\end{enumerate}
\end{Def}

\begin{Rem}
For the classification of one-dimensional Nash groups,  see \cite{MS}.
\end{Rem}

Nash submanifolds and subgroups are defined in obvious manner:

\bd\Label{nash}
A {\em Nash submanifold} in a Nash manifold is any real-analytic
submanifold, which can be defined locally in 
a neighborhood of each its point by 
Nash functions $f_1=\ldots=f_m=0$
satisfying $df_1\wedge\ldots\wedge df_m\ne 0$.
A {\em Nash subgroup} $H$ of a Nash group is any subset
which is both Nash submanifold and a subgroup.
\ed

In particular, Nash group is always a Lie group
with finitely many connected components
and a Nash subgroup of a Nash group is always a real-analytic Lie subgroup.

Given a Nash manifold $M$, we call a subset
$S\subset M$  semi-algebraic if  it  has semi-algebraic intersection with every 
Nash coordinate chart.

\bl\Label{nash-sub}
Let $H$ be a subgroup in a Nash group $G$.
Assume that $H$ is also a semi-algebraic subset of $G$.
Then $H$ is a Nash subgroup of $G$.
\el

\bpf
Consider a Nash coordinate chart $\phi\colon U\to V\subset\R^n$ in $G$ around a point $g_0\in H$. Since $H$ is semi-algebraic, $S:=\phi(H\cap U)$ is semi-algebraic in $\R^n$. By Proposition~\ref{strat}, $S$ admits a finite stratification into disjoint semi-algebraic sets which are real-analytic submanifolds. Let $A\subset S$ be one of these submanifolds of the highest dimension and choose a point $a\in A$. Then 
$H$ is a Nash submanifold of $G$ in a neighborhood of $a':=\phi^{-1}(a)$.
Since $H$ is subgroup of $G$, it is also a Nash submanifold in a neighborhood of its every point. Hence $H$ is a Nash subgroup.
\epf

\bc\Label{proj-group}
Let $\phi\colon G\to G'$ be a Nash homomorphism between two Nash groups.
Then for every Nash subgroup $H\subset G$, the image $\phi(H)$ is a Nash subgroup of $G'$.
\ec

\bpf
By Theorem~\ref{TS} applied in coordinate charts,
$\phi(H)$ is a finite union of semi-algebraic subsets of $G'$
and is therefore semi-algebraic. Then $\phi(H)$
is a Nash subgroup of $G'$ by Lemma~\ref{nash-sub}.
\epf

\bl\Label{stabil-group}
Let $(H_{m})$ be a descreasing sequence of Nash subgroups of a Nash group $G$, i.e.\ $H_{m}\supset H_{m+1}$. Then $(H_{m})$ stabilizes, i.e.\ $H_{m}=H_{m+1}$ for $m$ sufficiently large.
\el

\br
Note that a decreasing sequence of Lie subgroups need not stabilize, e.g.\ 
take $H_{m}:=m\Z$.
\er

\bpf[Proof of Lemma~\ref{stabil-group}]
Clearly the dimension $\dim H_m$ stabilizes after some $m=m_0$.
Since for $m\ge m_0$, every subgroup $H_m\subset H_{m_0}$
has the same dimension, it must coincide with a union of some
of the connected components of $H_{m_0}$.
Since $H_{m_0}$ has only finitely many connected components,
it has  only finitely many possible Nash sugroups of the same dimension.
Hence the sequence $H_m$ must terminate.
\epf

\section{The formal division algorithm}
The second main ingredient is the formal division algorithm.
We closely follow the article \cite{BM87} of Bierstone and Milman, where the reader is 
referred for further details.

Let $\K$ be a field and $\K[[t]]$ the ring of all power series in $t=(t_{1},\ldots,t_{n})$.
Consider the order on the set of monomials $ct^{\a}=ct_{1}^{\a_{1}}\ldots t_{n}^{\a_{n}}$ 
induced by the lexicographic order from the right on the $(n+1)$-tuples
$(\a_{1},\ldots,\a_{n},|\a|)$. 

\bd\Label{initial}
For every power series $f\ne0$,
its {\em initial exponent} is the multi-index $\a$ of the smallest nonzero monomial
in the expansion of $f$. 
Given an ideal  $I\subset\K[[t]]$,
its {\em diagram of initial exponents} $\5N(I)$
is the set of all initial exponents of all nonzero elements $f\in I$.
\ed

\bt[Grauert, Hironaka]\Label{div-alg}
Let $g^{1},\ldots,g^{k}\in\K[[t]]$ be nonzero elements
with initial exponents $\a^{1},\ldots,\a^{k}$ respectively
and let $f\in\K[[t]]$ be another element.
Then there exist $q^{1},\ldots,q^{k},r\in\K[[t]]$ such that
$$f=q^{1}g^{1}+\ldots+ q^{k}g^{k} + r$$
and $r$ has in its expansion no nonzero monomials $c_{\a}z^{\a}$
with 
$$\a\in \bigcup_{j=1}^{k} (\a^{j}+\N^{n}).$$
\et

Since $I$ is invariant under multiplication with monomials,
it follows that $\5N(I)=\5N(I)+\N^{n}$.
We have the following elementary lemma:

\bl\Label{terminate}
Any increasing sequence $\5N_1\subset\5N_2\subset \ldots$
of subsets $\5N_k\subset\N^n$ satisfying $\5N_k=\5N_k+\N^n$
terminates.
\el

\bpf
We proceed by induction on $n$. The statement is obvious for $n=1$.
Suppose it holds for $n=l$ and consider the case $n=l+1$.
Suppose by contradiction that $(\5N_k)$ is a strictly increasing sequence of subsets of $\N^{l+1}$ satisfying the assumptions of the lemma. 
For each $a\in\N$, set
$$\5N_k^{a}:=\{\a\in\N^l : (\a,a)\in\5N_k\}, \quad \5N^a:=\cup_s\5N_s^a.$$
Since $\5N_k=\5N_k+\N^{l+1}$, it follows that 
$\5N_k^a\subset \5N_k^{a+1}$ and hence $\5N^a\subset\5N^{a+1}$.
Then, by the induction assumption, the sequence $(\5N^a)$ 
terminates for some $a=a_0$. Using the induction assumption again for each $a=0,\ldots, a_0$, we conclude that there exists $k_0$
such that each sequence $(\5N_k^{a})_k$ terminates after $k=k_0$. 
That is, assuming $k>k_0$, we have $\5N_k^a=\5N_{k_0}^a=\5N^a$ for all $a=0,\ldots,a_0$. Furthermore, with the same assumption, for every $a>a_0$, 
we have
$$\5N_k^a\subset\5N^a=\5N^{a_0}=\5N_{k_0}^{a_0}\subset \5N_{k_0}^a$$
proving that the sequence $(\5N_k^a)_k$ terminates after $k=k_0$.
Thus $\5N_k^a=\5N_{k_0}^a$ holds for every $a$, hence $\5N_k=\5N_{k_0}$.
That is, the sequence $(\5N_k)$ terminates and the proof is complete.
\epf

\bc\Label{b}
Any subset $\5N\subset \N^n$ with $\5N=\5N+\N^n$
contains a finite subset $\5B\subset\5N$ with $\5N=\5B+\N^n$.
\ec

The minimal subset $\5B$ with that property is called the {\em set of vertices} in \cite[1.4]{BM87}.

\bpf[Proof of Corollary~\ref{b}]
Assume by contradiction, that for any finite subset $\5B\subset\5N$,
we have $\5N\ne\5B+\N^n$. 
Then we can construct iductively a sequence $(\b^k)$ in $\5N$
such that 
$$\b^{k+1}\notin \5N_k:=\{\b^1,\ldots,\b^k\}+\N^n$$
for every $k$.
Hence the sequence $(\5N_k)$ satisfy the assumptions of Lemma~\ref{terminate}
but does not terminate,
which is a contradiction.
\epf

Now given an ideal $I\subset\K[[t]]$, let $\5B=\{\b^{1},\ldots,\b^{k}\}\subset\5N(I)$ be any finite subset satisfying the conclusion of  Corollary~\ref{b} and choose  any $g_{1},\ldots,g_{k}\in I$ whose initial exponents are $\b^{1},\ldots,\b^{k}$ respectively. 
Then Theorem~\ref{div-alg} yields:

\bc\Label{cor-division}
Let $I\subset\K[[t]]$ be an ideal.
Then for every $f\in\K[[t]]$, there exist $g\in I$ and $r\in\K[[t]]$
such that $f=g+r$ and $r$ has in its expansion no nonzero monomials $c_{\a}t^{\a}$
with $\a\in\5N(I)$.
\ec

In the following we denote by $\5m\subset\K[[t]]$
the maximal ideal consisting of all formal power series vanishing at $0$.

\bp\Label{ideal-finite}
Let $I\subset\K[[t]]$ be an ideal and $f\in\K[[t]]$ a formal power series
with $f\in I+ \5m^{k}$ for every $k$. Then $f\in I$.
\ep

\bpf
In view of Corollary~\ref{cor-division}, we may assume
that $f$ has no monomials $c_\a t^\a$ in its expansion
with $\a\in\5N(I)$. Assume by contradiction that $f\ne0$.
Then choose any $k$ with $f\notin\5m^k$.
By the assumption, there exists $g\in I$ with $f-g\in\5m^k$.
Then $f$ and $g$ in their expansions have the same monomials 
of order less than $k$. In particular, $g$ in its expansion 
has a nonzero monomial of order less than $k$
but none of those monomials $c_\a t^\a$ satisfies $\a\in\5N(I)$.
By Definition~\ref{initial}, the initial exponent of $g$
does not belong to $\5N(I)$, which is a contradiction with the construction of $\5N(I)$, because $g\in I$.
\epf

Proposition~\ref{ideal-finite} can be restated as
$j^{k}f\in j^{k}I$ for all $k$ implies $f\in I$.
Here $j^{k}f$ is the $k$-jet of $f$ (at $0$),
which is the equivalence class of $f$,
where two formal power series are equivalent
if they coincide up to order $k$.
Furthermore, $j^{k}I:=\{j^{k}g:g\in I\}$.
We shall use the following consequence of Proposition~\ref{ideal-finite}:

\bc\Label{ideal-equiv}
Let $\Phi\colon(\K^{m},0)\to (\K^{n},0)$ be a formal map
and $I\subset\C[[t_{1},\ldots,t_{m}]]$, $I'\subset \C[[t_{1},\ldots,t_{n}]]$ be two ideals.
If $\Phi^{*}(j^{k}I')\subset j^{k}I$ for all $k$, then $\Phi^{*}I'\subset I$.
\ec

\bpf
Fix any $g\in I'$. Then the assumption $\Phi^{*}(j^{k}I')\subset j^{k}I$
implies that $g\circ\Phi\in I+\5m^k$ for every $k$.
By Proposition~\ref{ideal-finite}, we have $g\circ\Phi\in I$.
Since $g\in I'$ was arbitrary, we obtain the desired conclusion.
\epf

\section{Formal and finite order equivalences}
Here we give a proof of the main theorem 
describing the set of all formal equivalences
between two families of ideals in $\R[[x,y]]$,
where $z=x+iy\in\R^{n}+i\R^{n}$.
We continue using the notation 
$j^{k}f$ and $j^{k}I$ from the previous section
for $k$-jets (at $0$) of  a formal power series $f$ and an ideal $I$ respectively.
(Note that $k$-jets here only make sense at $0$.)
We further denote by $\6G^{k}$
the group of all invertible $k$-jets of formal maps
$\Phi\colon (\C^{n},0)\to (\C^{n},0)$.
It is easy to see that $\6G^{k}$
has a natural structure of a Nash group (see \S\ref{Nash} for this notion).
Given any $k$-jet $\L\in \6G^{k}$, its $l$-jet $j^{l}\L\in \6G^{k}$ for $l<k$
is defined in an obvious way by truncation.

\bt\Label{main}
Let $(I_{\a})_{\a\in A}$, $(I'_{\a})_{\a\in A}$ be two given families of ideals in $\R[[x,y]]$
that are equivalent of any finite order. Then there exists a sequence of Nash subgroups
$H^{k}\subset  \6G^{k}(\C^{n})$, $k=1,2,\ldots$,
and of right $H^{k}$-cosets $R^{k}\subset  \6G^{k}(\C^{n})$
such that the following hold:
\begin{enumerate}
\item[(i)] $j^{k}H^{k+1}=H^{k}$ and $j^{k}R^{k+1}=R^{k}$ for all $k$;
\item[(ii)] a formal map $\Phi\colon (\C^n,0)\to(\C^n,0)$ is a formal equivalence between $(I_{\a})$ and $(I'_{\a})$
if and only if $j^{k}\Phi\in R^{k}$ for all $k$;
\item[(iii)] for every $k$ there exists $l$ such that if $\Phi$ is an equivalence
of order $l$ between $(I_{\a})$ and $(I'_{\a})$, then
there exists a formal equivalence $\2\Phi$ between $(I_{\a})$ and $(I'_{\a})$ with $j^{k}\2\Phi=j^{k}\Phi$.
\end{enumerate}
\et

\bpf[Proof of Theorem~\ref{main}]
For every integer $k\ge 1$, consider $k$-jets 
$j^{k}I_{\a},j^{k}I'_{\a}\subset\R[[x,y]]/\5m^{k}$, $\a\in A$.
Let $G^{k}\subset \6G^{k}(\C^{n})$
be the subgroup of all $k$-jets preserving each $j^{k}I'_{\a}$, $\a\in A$.
Then $G^{k}$ is a real algebraic subgroup of $\6G^{k}(\C^{n})$.
Let $J^{k}\subset \6G^{k}(\C^{n})$ be the subset of all $k$-jets
sending $j^{k}I_{\a}$ onto $j^{k}I'_{\a}$ for each $\a\in A$.
Note that, by definition, $\Phi$ is an equivalence of order $k+1$
between the families $(I_{\a})$ and $(I'_{\a})$
 if and only if
$j^{k}\Phi\in J^{k}$.
According to our assumption, $J^{k}\ne\emptyset$ for every $k$.
Then $J^{k}$ is a right $G^{k}$-coset,
i.e.\ $J^{k}=\L\cdot G^{k}$, where $\L\in J^{k}$ is any element.

We next consider the truncation map $\6G^{l}(\C^{n})\to \6G^{k}(\C^{n})$, $\L\mapsto j^{k}\L$,
for $l\ge k$ and the image subgroups $j^{k}G^{l}\subset \6G^{k}(\C^{n})$.
By Corollary~\ref{proj-group}, each $j^{k}G^{l}$ is a Nash subgroup of $\6G^{k}(\C^{n})$.
Furthermore, if a $l$-jet $\L$ preserves each $j^{l}I_{\a}$, $\a\in A$,
its $k$-truncation $j^{k}\L$ preserves each $j^{k}I_{\a}$.
Hence $j^{k}G^{l}\supset j^{k}G^{l+1}$, i.e.\  $(j^{k}G^{l})_{l}$ is a decreasing sequence of
Nash subgroups of $\6G^{k}(\C^{n})$.
Similarly, $(j^{k}J^{l})_{l}$ is a decreasing sequence of subsets.
By Lemma~\ref{stabil-group}, there exists $l=l(k)\ge k$
such that 
\begin{equation}\Label{stab-g}
j^{k}G^{l} = j^{k}G^{l'} \quad\text{ for all }\quad l'\ge l.
\end{equation}
Since the right $G^{l}$-action on $\6G^{l}(\C^{n})$ commutes with truncation,
every $j^{k}J^{l}\subset \6G^{k}(\C^{n})$ is a right $j^{k}G^{l}$-coset.
Since the sequence $(j^{k}J^{l})_{l}$ is decreasing, \eqref{stab-g}
implies 
\begin{equation}\Label{stab-j}
j^{k}J^{l} = j^{k}J^{l'} \quad\text{ for all }\quad l'\ge l.
\end{equation}

We now set
\[
H^{k}:= j^{k}G^{l(k)} =  \bigcap_{l\ge k}j^{k}G^{l} \subset \6G^{k}(\C^{n}), \quad
R^{k}:=  j^{k}J^{l(k)} = \bigcap_{l\ge k}j^{k}J^{l}\subset \6G^{k}(\C^{n}).
\]
Then each $H^{k}$ is a Nash subgroup
and each $R^{k}$ is a right $H^{k}$-coset.
We claim that $H^k$ and $R^k$ satisfy the conclusions of the Theorem.

Indeed, choosing $l=\max(l(k),l(k+1))\ge k+1$, 
we have $H^{k}=j^{k}G^{l}$ and $H^{k+1}=j^{k+1}G^{l}$ as consequence of \eqref{stab-g}.
Therefore, 
$$j^{k}H^{k+1}=j^{k}(j^{k+1}G^{l})=j^{k}G^{l}=H^{k}$$
proving the first identity in (i).
The proof of the second identity is completely analogous.

Now let $\Phi\colon (\C^{n},0)\to (\C^{n},0)$ be a formal equivalence
between $(I_{\a})$ and $(I'_{\a})$. 
Then we have $j^{k}\Phi\in J^{k}$ by the construction of $J^k$ and therefore 
$$j^{k}\Phi=j^{k}(j^{l(k)}\Phi)\in j^{k}J^{l(k)}=R^{k}.$$
Vice versa, let $\Phi$ be a formal map satisfying $j^{k}\Phi\in R^{k}$
for all $k$. Since $R^{k}\subset J^{k}$, we also have $j^{k}\Phi\in J^{k}$,
i.e.\ $\Phi$ is an equivalence of order $k$ for every $k$.
The latter property means that $\Phi^*(j^k I'_\a)=j^k I_\a$
for every $k$ and $\a$. Now Corollary~\ref{ideal-equiv}
applied for $\Phi$ and $\Phi^{-1}$,
implies that $\Phi^* I'_\a = I_\a$ for all $\a$,
i.e.\ $\Phi$ is a formal equivalence
between the families $(I_\a)$ and $(I'_\a)$.
This proves (ii).

Finally to show (iii), for every $k$, choose $l=l(k)$ as above
with the property that $R^k=j^k J^{l}$
and let $\Phi$ be an equivalence between $(I_{\a})$ and $(I'_{\a})$ of order $l+1$.
Then by our construction, we have $j^l\Phi\in J^l$
and hence $j^k\Phi\in R^k$.
In view of (i), we can construct inductively a sequence
$\L_m\in R^m$, $m\ge k$, satisfying $\L_k=j^k\Phi$
and $\L_m=j^m\L_{m+1}$.
Then the sequence $(\L_m)_{m\ge k}$ determines a unique formal map 
$\2\Phi\colon (\C^{n},0)\to (\C^{n},0)$ with $j^m\2\Phi = \L_m$
for all $m\ge k$. In particular, $\2\Phi$
satisfies $j^k\2\Phi=j^k\Phi$ and $j^m\2\Phi\in R^m$ for all $m$.
In view of (ii), the latter property implies that $\Phi$
is a formal equivalence
between $(I_{\a})$ and $(I'_{\a})$ as desired.
\epf

\section{Equivalences for self-maps and vector fields}\label{fields}
We next consider the set $\6E_{n}$ of all germs at $0$ of real-analytic self-maps 
$F$ of $\C^{n}$ preserving $0$
and equivalence relations on $\6E_{n}$ by conjugations.
More precisely, let $\Phi\colon (\C^{n},0)\to(\C^{n},0)$ be as before a formal invertible map given by 
a power series in $z=(z_{1},\ldots,z_{n})\in\C^{n}$
vanishing at $0$.

\begin{enumerate}
\item
$\Phi$ is said to be a {\em formal equivalence} between two germs
$F,G\in\6E_{n}$ if $G=\Phi\circ F\circ \Phi^{-1}$.
\item $\Phi$ is said to be an {\em equivalence of order $k$}  between two germs
$F,G\in\6E_{n}$ if $G-\Phi\circ F\circ \Phi^{-1}$ vanishes of order at least $k$ at $0$, we write
$G\sim_{k}\Phi\circ F\circ \Phi^{-1}$.
\item $\Phi$ is said to be a {\em formal equivalence} (resp.\ {\em equivalence of order $k$}) between two families
$(F_\a)_{\a\in A}$ and $(G_\a)_{\a\in A}$ of germs in $\6E_{n}$ if 
$G_{\a}=\Phi\circ F_{\a}\circ \Phi^{-1}$  (resp.\ $G_{\a}\sim_{k}\Phi\circ F_{\a}\circ \Phi^{-1}$)  for all $\a\in A$.
\end{enumerate}

We have the following analogue of Theorem~\ref{maincor}: 
\bt\Label{maincor1}
Two families of germs of real-analytic self-maps of $\C^{n}$ preserving $0$ are formally equivalent
if and only if they are equivalent of any finite order.
\et

Theorem~\ref{maincor1} is obtained as a direct consequence of the following analogue of Theorem~\ref{main} for germs of real-analytic self-maps:
\bt\Label{main1}
Let $(F_{\a})_{\a\in A}$, $(F'_{\a})_{\a\in A}$ be two given families of 
germs at $0$ of real-analytic self-maps of $\C^{n}$ fixing $0$
that are equivalent of any finite order. Then there exists a sequence of Nash subgroups
$H^{k}\subset  \6G^{k}(\C^{n})$, $k=1,2,\ldots$,
and of right $H^{k}$-cosets $R^{k}\subset  \6G^{k}(\C^{n})$
such that the following hold:
\begin{enumerate}
\item[(i)] $j^{k}H^{k+1}=H^{k}$ and $j^{k}R^{k+1}=R^{k}$ for all $k$;
\item[(ii)] a formal map $\Phi\colon (\C^n,0)\to(\C^n,0)$ is a formal equivalence between $(F_{\a})$ and $(F'_{\a})$
if and only if $j^{k}\Phi\in R^{k}$ for all $k$;
\item[(iii)] for every $k$ there exists $l$ such that if $\Phi$ is an equivalence
of order $l$ between $(F_{\a})$ and $(F'_{\a})$, then
there exists a formal equivalence $\2\Phi$ between $(F_{\a})$ and $(F'_{\a})$ with $j^{k}\2\Phi=j^{k}\Phi$.
\end{enumerate}
\et

\bpf
The proof follows closely the line of the proof of Theorem~\ref{main}.
We write $j^{k}\6E_{n}$ for the space of all $k$-jets (at $0$) of elements in $\6E_{n}$.
Then consider $j^{k}F_{\a}, j^{k}F'_{\a}\in j^{k}\6E_{n}$ and
let $G^{k}\subset\6G^{k}$ be the subgroup of all $k$-jets $\L$
preserving each $j^{k}F_{\a}$, i.e.\ such that $\L\circ j^{k}F_{\a}\circ\L^{-1}=j^{k}F_{\a}$ for all $\a$.
Let further $J^{k}\subset\6G^{k}$ be the subset of all $k$-jets $\L$
sending $j^{k}F_{\a}$ to $j^{k}F'_{\a}$, i.e.\ such that $\L\circ j^{k}F_{\a}\circ\L^{-1}=j^{k}F'_{\a}$ for all $\a$.
As in the proof of Theorem~\ref{main}, we note that $\Phi$ is an equivalence of order $k+1$
between $(F_{\a})$ and $(F'_{\a})$ if and only if $j^{k}\Phi\in J^{k}$.
Hence $J^{k}\ne\emptyset$ for all $k$ and $J^{k}$ is a right $G^{k}$-coset.

Again, each $G^{k}$ is a Nash subgroup of $\6G^{k}$.
The rest of the proof repeats that of the proof of Theorem~\ref{main}.
\epf

A closely related important situation is that of equivalences for {\em vector fields}.
Here we regard $\Phi$ as a change or coordinates,
so that $\Phi$ is an equivalence between two vector fields $\xi$ and $\xi'$
if $\xi'=\Phi_{*}(\xi\circ\Phi^{-1})$. The notions of formal and finite order equivalences
between real-analytic vector fields
are defined in an obvious fashion analogously to the case of real-analytic self-maps above.
Now the study of equivalences for vector fields can be reduced to the
case of equivalences for self-maps by considering their flows.

Recall that the {\em flow} of a vector field $\xi$ in $\C^{n}$
is a smooth one-parameter family of local self-maps $F_{t}$ of $\C^{n}$
such that $F_{0}=\id$ and $\frac{d}{dt} F_{t}=\xi\circ F_{t}$.
Since the equivalence problem for non-singular germs of vector fields
(i.e.\ those not vanishing at the reference point) is trivial,
we restrict here only to germs at $0$ of singular real-analytic vector fields $\xi$,
i.e.\ such that $\xi(0)=0$.
Then the local flow $F_{t}$ is defined for all $t$ and consists of germs $F_{t}\in\6E_{n}$.
Furthermore, by expanding $F_{t}$ and $\xi$ into power series in $z$,
it is easy to see that two germs of (real-analytic) vector fields 
are formally equivalent (resp.\ equivalent of finite order)
if and only if their flows are formally equivalent (resp.\ equivalent of finite order).
Hence we obtain the following direct corollary of Theorem~\ref{maincor1}:

\bc\Label{maincor2}
Two families of germs of singular real-analytic vector fields in $\C^{n}$ are formally equivalent
if and only if they are equivalent of any finite order.
\ec

\section{Two sets of curves that are not formally equivalent
but equivalent of any finite order}\Label{example}
Our discussion here over complex numbers
can be repeated word for word for real numbers
without any change.

We shall consider complex plane curves in $\C^2$ passing through $0$ and
given in the coordinates $(z,w)$ by $w=\phi(z)$,
where $\phi$ is a polynomial.
Each of the two sets of such curves
will be indexed by two integers $(m,n)\in\N\times\Z$.
That is, we define two families 
$\{w=\phi_{m,n}(z)\}$ and $\{w=\psi_{m,n}(z)\}$.

We first construct inductively a sequence of integers $c_m$, $m\ge 1$,
such that the subsets $S_m:=2^m\Z+c_m\subset\Z$ satisfy
\begin{equation}\Label{cond}
S_m\supset S_{m+1}, \quad 
\bigcap_{m\ge 1} S_m = \emptyset.
\end{equation}
We first put $c_1:=1$, so that $S_1$ is the set of all odd integers.
Suppose that we have already constructed $c_1,\ldots,c_k$
with $S_1\supset\ldots\supset S_k$.
Then $0\notin S_k$ and hence $S_k$ has the maximum negative element $a_k<0$
and the minimum positive element $b_k>0$ such that $b_k-a_k=2^k$.
We put $c_{k+1}:=a_k$ if $|a_k|>b_k$ and $c_{k+1}:=b_k$ otherwise.
Obviously $S_k\supset S_{k+1}$. Furthermore, we have $2^{k-1}\le |c_{k+1}|< 2^k$. Then it follows that $|l|\ge 2^{k-1}$ for any $l\in S_{k+1}$.
Using the above procedure, we construct $c_m$ inductively
satisfying the first condition in \eqref{cond}
and such that $|l|\ge 2^{m-2}$ for any $l\in S_{m}$.
The latter property immediately implies
the second condition in \eqref{cond} as desired.

We now set
$$\phi_{m,n}(z):=2^m nz+z^{m+1}, 
\quad \psi_{m,n}(z):=(2^m n+c_m)z+z^{m+1}.$$
Our main conclusion of this section is the following:

\bp
The sets of germs at $0$ of the curves $\{w=\phi_{m,n}(z)\}$
and $\{w=\psi_{m,n}(z)\}$
are equivalent of any finite order
but not formally equivalent.
\ep

\bpf
As a consequence of the first inclusion property in \eqref{cond},
it is easy to see that for each $k$, the map $\Phi\colon (z,w)\mapsto (z,w+c_kz)$
defines an equivalence between the two sets of order $k+2$.
Indeed, since $c_k\in S_m$ for each $m\le k$, the map $l\mapsto l+c_k$
defines a bijection between $2^m\Z$ and $2^m\Z+c_k=S_m$.
Therefore, for $m\le k$, $\Phi$ maps any curve $w=\phi_{m,n}(z)$
into $w=\psi_{m,n'}(z)$ for suitable $n'$ and $\Phi^{-1}$ maps
any curve $w=\psi_{m,n}(z)$ into $w=\phi_{m,n'}(z)$ for suitable $n'$.

On the other hand, if $m>k$, we have $\phi_{m,n}(z)=2^mnz+O(|z|^{k+2})$
and hence $\phi_{m,n}(z)=\phi_{k,n'}(z)+O(|z|^{k+2})$
for suitable $n'$. Therefore $\Phi$ maps $w=\phi_{m,n}(z)$
into $w=\psi_{k,n''}(z)$ up to order $k+2$ for suitable $n''$.
Vice versa, it follows from the inclusion property in \eqref{cond}
that for $m>k$, $\psi_{m,n}(z)=\psi_{k,n'}(z)+O(|z|^{k+2})$ for some $n'$.
Hence $\Phi^{-1}$ maps $w=\psi_{m,n}(z)$ into $w=\phi_{k,n''}(z)$
for suitable $n''$.

We now show that the two sets are not formally equivalent.
By contradiction, suppose that we have 
a formal invertible map sending each $w=\phi_{m,n}(z)$
into some $w=\psi_{m',n'}(z)$.
Consider the curves $w=\phi_{m,0}(z)$,
whose tangent spaces at $0$ are all equal to $\{w=0\}$.
Therefore, the second set of curves $w=\psi_{m,n}(z)$
must contain an infinite collection of curves whose tangent
spaces at $0$ all coincide. However, it follows from 
the second condition in \eqref{cond}
that the latter is impossible.
\epf

\end{document}